\input amstex
\magnification=1200
\loadmsam
\loadmsbm
\loadeufm
\loadeusm
\UseAMSsymbols

\hsize=6.25truein
\hoffset=0.15truein
\vsize=9truein
\voffset=-0.2truein

\def\leftitem#1{\item{\hbox to\parindent{\enspace#1\hfill}}}

\def\boxit#1#2{\hbox{\vrule
	\vtop{%
	\vbox{\hrule\kern#1%
	\hbox{\kern#1#2\kern#1}}%
	\kern#1\hrule}%
	\vrule}}

\def\leaderfill{\leaders\hbox to 1em{\hss.\hss}\hfill}

\parskip=\medskipamount
\document

\input epsf

\centerline{\bf A Presentation of the Mapping Class Groups}
\bigskip
\centerline{Feng Luo}

\centerline{\it Dept. of Math., Rutgers University,  New Brunswick, NJ 08903 \rm}
\centerline{e-mail: fluo\@math.rutgers.edu}

{\bf Abstract.} Using the works of Gervais, Harer, Hatcher and Thurston  and
others, we show that the mapping class
group of  a compact orientable surface has a presentation 
so that the generators are the set of all Dehn twists and 
the relations are supported in subsurfaces homeomorphic to 
the one-holed torus or the four-holed sphere. It turns out that
all the relations were discovered by Dehn in 1938.

\noindent
\S 1. {\bf Introduction}

Let $\Sigma$ = $\Sigma_{g,n}$ be a compact oriented surface of genus $g$ with
$r$ boundary components and $\Cal M_{g,r} = \Cal M(\Sigma)$ be the
mapping class group Home($\Sigma, \partial \Sigma)/$Iso where homeomorphisms
and isotopies leave points on $\partial \Sigma$ fixed. A presentation of
the mapping class groups in terms of Dehn twists was obtained in the
fundamental paper by Hatcher and Thurston [HT] where relations are
supported in subsurfaces homeomorphic to $\Sigma_{2,3}$. 
Using the work of [HT], [Ha], and [Wa], Gervais in [Ge] obtains 
a presentation of $\Cal M_{g,r}$ where relations are supported in 
subsurfaces homeomorphic to $\Sigma_{1,2}$. 
The goal of this note is to simplify Gervais' presentation and show 
that relations are actually supported in  subsurfaces 
homeomorphic to $\Sigma_{1,1}$ or $\Sigma_{0,4}$.

We begin by introduce some notations. Let $\Cal S$ = $\Cal S(\Sigma)$ be the
set of isotopy classes of simple loops on the surface $\Sigma$. Given $\alpha$,
$\beta$ in $\Cal S$, define $I(\alpha, \beta) = $ min$\{|a \cap b|:
a \in \alpha, b \in \beta\}$. We use $\alpha \cap \beta = \emptyset$ to
denote $I(\alpha, \beta) = 0$; use $\alpha \perp \beta$ to denote
$I(\alpha, \beta) =1$; and use $\alpha \perp_0 \beta$ to denote
$I(\alpha, \beta) = 2$ so that their algebraic intersection number is
zero. If $a, b$ are two arcs intersecting  transversely at a point $p$,
then the \it resolution of \rm  $a \cup b$  \it at $p$ from $a$ to $b$ \rm
is defined as  follows. Fix any orientation  on $a$ and use the orientation
on the surface to determine an orientation on $b$. Then resolve
the intersection according to the orientations (see figure 1). 
The resolution is
evidently independent of the choice of the orientations on $a$. If
$\alpha \perp \beta$ or $\alpha \perp_0 \beta$, take $a \in \alpha$,
$b \in \beta$ so that $|a \cap b| = I(\alpha , \beta)$. Then the
curve obtained by resolving all intersection points in
$a \cap b$ from $a$ to $b$ is again a simple loop denoted by $ab$.
We define $\alpha \beta$ to be the isotopy class of $ab$. It follows
from the definition that when $\alpha \perp \beta$ then $\alpha \beta
\perp \alpha, \beta$ and when $\alpha \perp_0 \beta$ then $\alpha
\beta \perp_0 \alpha, \beta$. Also the Dehn twist along $\alpha$ applied
to $\beta$ can be recovered from the operation as follows: 
if  $\alpha \perp \beta$, then $D_{\alpha}(\beta) = \alpha \beta$;
and if $\alpha \perp_0 \beta$, then 
$D_{\alpha}(\beta) = \alpha (\alpha \beta)$.
Let $N(a)$ and $N(b)$ be two small regular neighborhoods of $a$ and $b$. Then
$N(a \cup b) = N(a) \cup N(b)$ is homeomorphic to $\Sigma_{1,1}$ when
$\alpha \perp \beta$ and to $\Sigma_{0,4}$ when $\alpha \perp_0 \beta$.
Let $\partial(\alpha, \beta)$ be the isotopy class of the curve system
$\partial N(a \cup b)$.

In terms of these notations, the result of Gervais is as follows.

\noindent
{\bf Theorem} (Gervais). \it For a compact oriented surface $\Sigma$, the
mapping class group $\Cal M$$(\Sigma)$ has the following presentation:

generators: \rm \{$D_{\alpha} : \alpha \in \Cal S(\Sigma)$\}.

\it relations: \rm \text{(I)} $D_{\alpha} = 1$ \it if $\alpha$ is the isotopy
class of the null homotopic loop. \rm

\text{(II)} $D_{\alpha} D_{\beta} = D_{\beta} D_{\alpha}$ \it
if $\alpha \cap \beta = \emptyset$. \rm

\text{(III)} $D_{\alpha \beta} = D_{\alpha}D_{\beta} D_{\alpha}^{-1}$  \it if
$\alpha \perp \beta$. \rm

\text{(III')}  $D_{D_{\alpha }(\beta)} =  D_{\alpha}D_{\beta} D_{\alpha}^{-1}$ 
 \it if $\alpha \perp_0 \beta$.\rm

\text{(IV)} $D_{\alpha}D_{\beta}D_{\alpha \beta} = D_{\partial(\alpha, \beta)}$
 \it if $\alpha \perp_0 \beta$. \rm

\text{(V')} $(D_{\alpha}D_{\beta}D_{\gamma})^4 = D_{\epsilon_1} D_{\epsilon_2}$
\it where $\alpha$, $\beta$, $\gamma$ and $\epsilon_i$ are as shown in 
figure 1. \rm

\midspace{0.1cm}
\centerline{\epsfbox{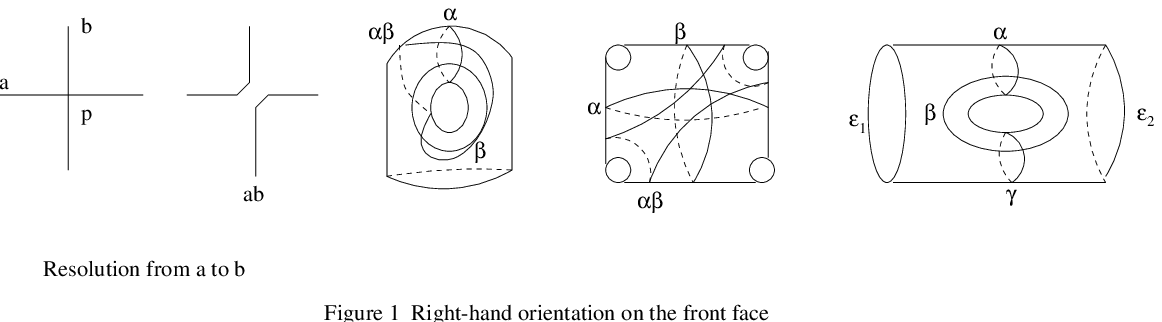}}
\midspace{0.1cm}

If  we set $\epsilon_1$ in (V') to be the trivial class, then (V') becomes
the relation
$$ (D_{\alpha} D_{\beta} D_{\alpha})^4 = D_{\partial(\alpha, \beta)}
\quad \quad \text{ if $\alpha \perp \beta.$} \tag V$$

Our observation is that relations (I), (II), (III), (IV), and (V)
form a complete set of relations.

\noindent
{\bf Theorem.} \it For compact oriented surface $\Sigma$, the mapping class
group $\Cal M$($\Sigma$) has a presentation where the generators are the
set of all Dehn twists and the relations are \rm  (I), (II), (III), (IV) \it
and \rm (V).

\noindent
\it Remarks. \rm 
1.  The relations (I)-(V) are well known.  Evidently relations 
(I), (II) hold.  Relation (IV)  is the
lantern relation which was discovered by Dehn ([De], pp333) in 1938 
and rediscovered independently by Johnson [Jo] in 1979. Relations
(III) and (V) were also discovered by Dehn ([De], pp287, and the last
sentence on pp310).

2. For surface $\Sigma$ $\cong \Sigma_{1,1}$ or $\Sigma_{0,4}$, let
$\Cal S'(\Sigma)$ be the set of isotopy classes of essential non-boundary
parallel simple loops on $\Sigma$. It is well known that there is a
bijection (the slope map) $\pi:$ $\Cal S'(\Sigma)$ $\to \hat \bold Q$ so that
 $\pi(\alpha) = p/q$,
$\pi(\beta) = p'/q'$ satisfy  $p'q - pq' = \pm 1$ if and only if
$\alpha \perp \beta$ or $\alpha \perp_0 \beta$. Futhermore, 
$\alpha \beta = (p+ \lambda p')/(q + \lambda q')$ where
$\lambda = pq' - p'q$. Thus the three classes $\alpha, \beta$
and $\alpha \beta$ in relations (III), (IV), and (V) form
an ideal triangle under the slope map in the modular configuration
$(\hat \bold Q, PSL(2, \bold Z))$.  This shows that the mapping 
class group can be reconstructed explicitly (in terms of
presentation) from the modular relation and the disjoint relation 
on $\Cal S(\Sigma)$.

\noindent
\S2. {\bf Proof of the Theorem}

As a convention, all surfaces drawn in the note have the right-hand
orientation in the front face.

It suffices to show that relations (III') and (V') are the
consequences of the relations (I)-(V).

To derive (III') that $D_{\gamma} = D_{\alpha}D_{\beta}D_{\alpha}^{-1}$
where $\gamma = D_{\alpha}(\beta)$ with $\alpha \perp_0 \beta$, we
note that $\alpha \beta \perp_0 \alpha$ and $\gamma = \alpha (\alpha \beta)$.
Furthermore, $\partial(\alpha, \beta) = \partial( \alpha, \alpha \beta)$.
Denote  $\partial(\alpha, \beta) $ by $\delta$. By  relation (IV) for $\alpha
\perp_0 \beta$ and $\alpha \beta \perp_0 \beta$, we obtain
$D_{\alpha}D_{\beta} D_{\alpha \beta} = D_{\delta}$ and
$D_{\alpha}D_{\alpha \beta} D_{\gamma} = D_{\delta}$. By relation (II),
$D_{\delta}$ commutes with $D_{\alpha}$, $D_{\beta}$, and $D_{\alpha \beta}$.
By cancelling $D_{\alpha \beta}$ and $D_{\delta}$ from the two
equations, we obtain $D_{\gamma} = D_{\alpha} D_{\beta} D_{\alpha}^{-1}$
which is the relation (III').

To derive the relation (V') from (I)-(V), we need two lemmas.

\noindent
{\bf Lemma 1.} \it If $\alpha \perp \beta$, then $(\alpha \beta) \alpha
= \alpha (\beta \alpha) = \beta$. In particular, as a consequnce of
relation \rm  (III), \it we obtain Artin's relation \rm
$$ D_{\alpha} D_{\beta} D_{\alpha} = D_{\beta} D_{\alpha} D_{\beta}
\quad \quad \text{if $\alpha \perp \beta$.} \tag VI$$

\noindent
PROOF (Well known).   To show  $(\alpha \beta) \alpha = \beta$, consider
figure 2.

\midspace{0.1cm}
\centerline{\epsfbox{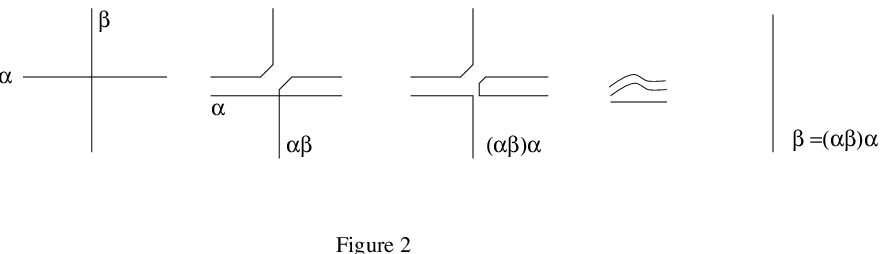}}
\midspace{0.1cm}

By (III) applied to $\alpha \perp \beta$, we obtain $D_{\alpha \beta}
= D_{\alpha}D_{\beta} D_{\alpha}^{-1}$.
Since $\alpha \beta \perp \alpha$ and  $(\alpha \beta) \alpha = \beta$, by
(III), we obtain $D_{\beta} = D_{\alpha \beta} D_{\alpha}
D_{\alpha \beta}^{-1}$. Combining these two equations, we obtain (VI).
$\square$

Now to show (V') that 
$(D_{\alpha}D_{\beta}D_{\gamma})^4 = D_{\epsilon_1} D_{\epsilon_2}$
where $\alpha$, $\beta$, $\gamma$ and $\epsilon_i$ are as in
figure 1,  we  need the following   lemma  (lemma 6.2, in [Lu])
whose proof is given in figure 3.

\noindent
{\bf Lemma 2.} \it Let $\delta$ be  $\partial (\alpha, \beta)$. Then
$\beta \gamma \perp  \delta \gamma$ and $(\beta \alpha) \alpha
= (\beta \gamma)(\delta \gamma)$. \rm

\midspace{0.1cm}
\centerline{\epsfbox{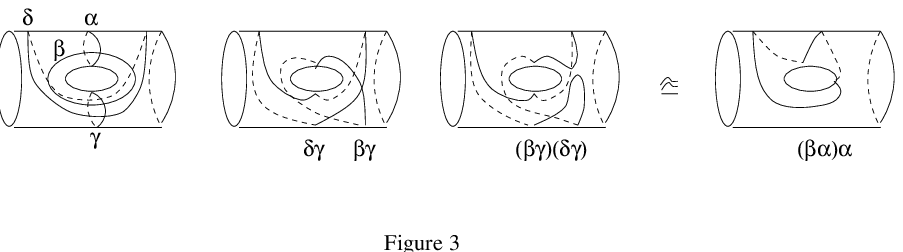}}
\midspace{0.1cm}

For simplicity, we use the  letters $ A,B,C,D,$ $E_i$ to denote
the Dehn twists on $\alpha, \beta, \gamma, \delta, \epsilon_i$ respectively.
For instance, (V') becomes $(ABC)^4 = E_1 E_2$.

By $\beta \perp \gamma$ and (III), we obtain
$$ D_{\beta \gamma} = BCB^{-1}. \tag 1$$
By $\alpha \perp \beta$, $\beta \alpha \perp \alpha$ and (III), we obtain
$D_{\beta \alpha} = BAB^{-1}$ and $D_{(\beta \alpha)\alpha}$$
= D_{\beta \alpha}D_{\alpha}D_{\beta \alpha}^{-1}$$ =$
$BAB^{-1}ABA^{-1}B^{-1}$. Using (VI) that $ABA = BAB$ and its
equivalent forms $BAB^{-1} = A^{-1}B A$, $ABA^{-1} = B^{-1}AB$, we obtain
$$ D_{(\beta \alpha)\alpha} = A^{-2}B A^2.  \tag 2$$
By $\alpha \perp \beta$ with $\partial(\alpha, \beta) = \delta$ and
$(V)$, we obtain
$$ D = (ABA)^4. \tag 3$$
By $\delta \perp_0 \gamma$ with $\partial (\delta, \gamma) = 
\epsilon_1 \cup \epsilon_2 \cup \alpha \cup \alpha$ and (IV), we obtain
$DCD_{\delta \gamma} = E_1 E_2 A^2$. By (II) that $E_i$ commutes
with $A,B,C,D$, we obtain
$$D_{\delta \gamma} = E_1E_2C^{-1}D^{-1}A^2.  \tag 4$$
Finally, by lemma 2 that $\beta \gamma \perp \delta \gamma$,
$(\beta \alpha) \alpha = (\beta \gamma)(\delta \gamma)$ and (III), 
we obtain $$D_{(\beta \alpha)\alpha} = D_{\beta \gamma}
D_{\delta \gamma}D_{\beta \gamma}^{-1}. \tag 5$$
Substitute (1)-(4) into (5) and use relation (II), we obtain
$$ A^{-2} B A^2 BCB^{-1} A^{-2}(ABA)^4 CBC^{-1} B^{-1} = E_1 E_2.  \tag 6$$

We claim that under relations (I)-(VI), the  left-hand side of
the equation (6) is $(BCA)^4$.  But $E_i$ commutes with $A$ by (II). 
Thus $(ABC)^4 = E_1 E_2$ after a conjugation by $A$.

Here is the calculation. In each  
step of the derivation below, we apply one of the relations (II) that
$AC=CA$, (VI) that $ABA =BAB$, $BCB=CBC$ or the cancellation 
law $XX^{-1} =1$ to the letters underlined.

$ A^{-2} B A^2 BCB^{-1} A^{-2}(ABA)^4 CBC^{-1} B^{-1}$

$ =A^{-2}BA^2 B C B^{-1} A^{-2} ABA ABA ABA \underline{ABA} \quad \underline{CBC^{-1}} B^{-1}$

$= A^{-2}BA^2 B C B^{-1} \underline{ A^{-2} A}BA ABA ABA  BA \underline{B B^{-1}} C \underline{B B^{-1}}$

$=  A^{-2}BA^2 B C  \underline{ B^{-1}A^{-1} B}A ABA ABA  BA C$

$= A^{-2}B A^2 BC AB^{-1} \underline{A^{-1} A}ABAABABAC$

$= A^{-1} \underline{ A^{-1} BA} ABCA  \underline{B^{-1} AB}AABABAC$

$ = \underline{A^{-1}BA}B^{-1}AB \underline{CAA}B \underline{A^{-1}A}A BABAC$

$=BAB^{-1}B^{-1} \underline{ABA}ACBABABAC$

$= BAB^{-1} \underline{B^{-1} B}ABACBABABAC$

$= BA \underline{B^{-1}AB}ACBABABAC$

$= BAAB \underline{A^{-1}A}CBABABAC$

$= BAA \underline{BCB}ABABAC$

$= BA \underline{AC}B \underline{CA}BABAC$

$= BAC \underline{ABA} CBA BAC$

$= BAC BA \underline{BCB}ABAC$

$= BAC BA CB \underline{C A}BAC$

$= (BAC)^4$

$= (BCA)^4$.

Thus $(BCA)^4 = E_1 E_2$. But $E_i$ commutes with $A$ by (II). Thus $(ABC)^4 
= E_1 E_2$ after a conjugation by $A$. This finishes the proof.

\centerline{\bf Reference}

[Bi] Birman, J.:  Braids, links, and mapping class groups.  Ann. of Math. Stud., 82, Princeton Univ. Press, Princeton, NJ, 1975.

[De] Dehn, M.: Papers on group theory and topology. J. Stillwell (eds.).
 Springer-Verlag, Berlin-New York, 1987.

[Ge] Gervais, S.: Presentation and central extensions of
mapping class groups. Trans. Amer. Math. Soc. 348 (1996), no. 8, 3097--3132.

[Har] Harer, J.: The second homology group of the mapping class group
of an orientable surface. Invent. Math. 72 (1983), 221-239.

[HT] Hatcher, A., Thurston, W.: A presentation for the mapping class group of a closed orientable surface. Topology 19 (1980), 221-237.

[Jo] Johnson, D.: Homeomorphisms of a surface which acts trivially on
homology. Proc. of Amer. Math. Soc. 75 (1979), 119-125.

[Li] Lickorish, R.: A representation of oriented combinatorial 3-manifolds. Ann.  Math.  72 (1962), 531-540.

[Lu] Luo, F.: Simple loops on surfaces and their interesection numbers, preprint.

[Wa] Wajnryb, B.: A simple presentation for the mapping class group of an orientable surface. Israel J. Math. 45 (1983), 157-174.

\end